\documentclass[a4paper,12pt,twoside]{amsart}
\usepackage{amssymb}

\usepackage{amsfonts,amssymb,amscd, amsmath}
\usepackage{graphicx}
 \usepackage[left=2.3cm,top=3.3cm,right=2.3cm,bottom=2.7cm]{geometry}
\usepackage{mathrsfs}

\usepackage[all,ps,cmtip]{xy}

\DeclareRobustCommand{\SkipTocEntry}[5]{}

\def\llra{\hbox to 10mm{\rightarrowfill}}

\def\lllra{\hbox to 15mm{\rightarrowfill}}

\def\cI{\mathcal{I}}

\def\cF{\mathcal{F}}
\def\cL{\mathcal{L}}
\def\cO{\mathcal{O}}

\def\cQ{\mathcal{Q}}

\let\tilde\widetilde

\DeclareMathOperator{\rank}{rank}

\DeclareMathOperator{\codim}{codim}
\DeclareMathOperator{\Pic}{Pic}

\newtheorem{lem}{Lemma}[section]
\newtheorem{theo}[lem]{Theorem}

\newtheorem{prop}[lem]{Proposition}

\theoremstyle{definition}

\newtheorem{rema}[lem]{Remark}

\theoremstyle{remark}
\newtheorem*{remark*}{Remark}
\newtheorem*{note*}{Note}

\begin{document}
\title{Syzygies of some rational homogeneous varieties}
\author{Zhi Jiang}
\address{Zhi Jiang, Shanghai Center for Mathematical Sciences, China}
\email{zhijiang@fudan.edu.cn}
\date{\today}
\maketitle
\begin{abstract}
In this paper,  we study syzygies of  rational homogeneous varieties. We extend Manivel's result that a $p$-th power of an ample line bundle on a flag variety satisfies Propery $(N_p)$  to many rational homogeneous varieties of type $B$, $C$, $D$, and $G_2$.
\end{abstract}
\section{Introduction}
In this paper we study the sygygies of rational homogeneous varieties. We assume that all varieties are over algebraically closed field of characteristic $0$.

Let $X$ be a smooth projective variety of dimension $n$ and let $L$ be a very ample line bundle. Let $R_L=\bigoplus_{m\geq 0} H^0(X, L^{\otimes m})$ be the section ring of $L$ and let $S_L=\bigoplus_{m\geq 0} S^mH^0(X, L)$ be the polynomial ring of $H^0(X, L)$. Then $R_L$ is naturally a $S_L$-module. By the Hilbert syzygy theorem, there exists a finite graded minimal free resolution $E^*$ of $R_L$:
$$\cdots \rightarrow E^1\rightarrow E^0\rightarrow R_L\rightarrow 0,$$
where $E^i=\bigoplus_{j} S_L(-a_{ij})$. The generators of $E^i$ are called (higher) syzygies of $L$.  Traditionally, if $E^0=S_L,$  $L$ is said to be projectively normal or normally generated and if $E^0=S_L$ and $E^1=\bigoplus_j S_L(-2)$, $L$ is said to be normally presented (see \cite{Mum}). Green and Lazarsfeld extended the terminology in \cite{GL}: if $E^0=S_L$ and $E^i=\bigoplus S_L(-i-1)$ for $1\leq i\leq p$,  one says that $L$ satisfies Property $(N_p)$.
In particular, $L$ satisfies Property $(N_0)$ means that it is projectively normal and $L$ satisfies Property $(N_1)$ iff it is normally presented.

 Consider the evaluation map $H^0(X, L)\otimes \cO_X\rightarrow L$ and denote by $M_L$ its kernel. It  is  a classical result that if $H^1(X, L^{\otimes k})=H^1(X, \wedge^iM_L\otimes L^{\otimes k})=0$ for all $k\geq 1$ and $1\leq i\leq p+1$, then $L$ satisfies property $(N_p)$ (see  \cite[Proposition 2.5]{AN}).

Despite the simple cohomological criterion for Property $(N_p)$,  it is in general a very difficult problem to compute the syzygies of a line bundle. There are satisfactory results about syzygies of ample line bundles on smooth projective curves (see \cite{G1} or  \cite[Section 1.8]{Laz}). In higher dimensions,  very few general results are known. In \cite{G2}, Green showed that $\cO_{\mathbb P^n}(d)$ satisfies Property $(N_p)$ for $p\leq d$. On the other hand, Ottaviani and Paoletti proved in \cite{OP} that for $n\geq 2$ and $d\geq 3$,  $\cO_{\mathbb P^n}(d)$ does not satisfy Property $(N_{3d-2})$. It is still an open question whether or not  $\cO_{\mathbb P^n}(d)$  satisfies Property $(N_{3d-3})$.

 The following theorem of Manivel proved in \cite{Man} is a generalization of Green's result about syzygies of projective spaces.

\begin{theo} Let $X$ be a   flag variety. Assume that $L$ is an ample line bundle on $X$, which can be expressed as $H^{\otimes l}\otimes M$, where $H$ is ample and $M$ is nef. Then $L$ satisfies Property $(N_p)$ when $l\geq p$.
\end{theo}

Manivel's result is the first step to answer a question of Fulton (see \cite[Problem 4.5]{EL}) asking to compute the syzygies of ample line bundles on rational homogeneous varieties via representation theory.

We continue to study Fulton's question in this article. We extend partially Manivel's result to other type rational homogeneous varieties. One of the main result is the following, which is a direct corollary of Theorem \ref{symplectic}.

 \begin{theo}
Let $X$ be a rational homogeneous variety of type $C_n$.
Let $L$ be a line bundle on $X$, which can be expressed as $H^{\otimes l}\otimes M$, where $H$ is an ample line bundle and $M$ is a nef line bundle. Then $L$ satisfies Property $(N_p)$ in the following cases:
\begin{itemize}
\item[(1)]  $\mathrm{rank} \Pic(X)\leq 2$ and $l\geq p$;
\item[(2)] or $l\geq \mathrm{max}\{p, \frac{p+1}{n}+\frac{n-3}{2}\}$.
\end{itemize}
 \end{theo}

Remark that  $(2)$ implies that for $p\geq \frac{n}{2}-1$, $L$ satisfies Property $(N_p)$ when $l\geq p$.

We also have  parallel statements with slight modifications for rational homogeneous varieties of type $B_n$, $D_n$, or $G_2$ (see Section 5).

 \noindent\textbf{Acknowledgements}.

The author is supported by the National Key Research and Development Program of China (No.~2020YFA0713200),  by NSFC for Innovative Research Groups
(No. 12121001), by the natural Science Foundation of China  (No.~11871155 and  No.~11731004), and by the Natural Science Foundation of Shanghai (No.~21ZR1404500).

\section{Flag varieties}
\subsection{The Borel-Weil-Bott theorem}

The Borel-Weil-Bott theorem is an important tool to study rational homogeneous varieties.
We will apply a special version of it. Fix $V$ a vector space of dimension $n$ and positive integers $n>n_1>\cdots >n_k>0$.
Let $ F:=\mathrm{Fl}(n_k, \ldots, n_1;V)$ be the flag variety parametrizing  linear subspaces $$V_k\subset V_{k-1}  \subset \cdots \subset V_1$$ contained in   $V$, where $\dim V_i=n_i$. We denote respectively $\Sigma_k\subset \Sigma_{k-1}  \subset \cdots \subset \Sigma_1$   the tautological sub-bundles  on $F$. By convention, $\Sigma_{k+1}=0$ and $\Sigma_0=V\otimes \cO_F$. Let $\cQ_i:=\Sigma_{i-1}/\Sigma_i$ for $1\leq i\leq k+1$.

 We denote by $L_i$ the determinant of $\cQ_i$ for $1\leq i\leq k$. We know that $$\Pic(F)=\mathbb Z L_1\oplus\cdots \oplus \mathbb ZL_k.$$ Given a line bundle $L=L_1^{\otimes a_1}\otimes \cdots \otimes L_k^{\otimes a_k}$,
 we know that $L$ is nef or basepoint-free iff $a_1\geq a_2\geq \cdots \geq a_k\geq 0$  and $L$ is ample or very ample iff $a_1> a_2> \cdots > a_k> 0$.

Let $r_i$ be the rank of $\cQ_i$ for $1\leq i\leq k+1$. Then $n=r_1+r_2+\cdots +r_{k+1}$. For a partition $\alpha_i=(a_{i1},\ldots, a_{i r_i})$ consisting of $r_i$ non-increasing sequence of integers, we denote by $S^{\alpha_i}(\cQ_i)$ the corresponding Schur power of $\cQ_i$. For instance, $S^{(k,0,\ldots,0)}\cQ_i=S^k\cQ_i$ is the $k$-th symmetric power of $\cQ_i$. We denote by $I_i=(1,\ldots,1)$ the partition of length $r_i$ for $1\leq i\leq k+1$ and $S^{I_i}\cQ_i=\wedge^{r_i}\cQ_i=L_i$.

 Let $\alpha=(\alpha_1, \cdots, \alpha_{k+1})$ be a vector with $n$ entries. We re-order the entries of $\alpha':=\alpha-(1,2, \ldots, n)$ to make it a non-increasing sequence $\alpha''$ and denote $\overline{\alpha}= \alpha''+(1, 2, \ldots, n)$. When $\overline{\alpha}$ is still non-increasing, $\alpha''$ is strictly decreasing and
   we denote by $$i(\alpha)=\#\{(i, j)\mid 1\leq i<j\leq n \;\text{and}\; \alpha_i-i<\alpha_j-j\}$$ the number of inversions of the sequence $\alpha'$.
The Borel-Weil-Bott theorem that we will apply frequently is  the following one (see \cite[Section 4.1]{W} for the general version).
\begin{theo}\label{BBW}Under the above assumptions,
\begin{itemize}
\item If $\alpha' $ contains two identical entries,   $H^i(F, S^{\alpha_1}\cQ_1\otimes \cdots \otimes S^{\alpha_{k+1}}\cQ_{k+1})=0$ for each $i\geq 0$;
\item If $\overline{\alpha}$ is non-increasing, then $H^i(F, S^{\alpha_1}\cQ_1\otimes \cdots \otimes S^{\alpha_{k+1}}\cQ_{k+1})=0$ for   $i\neq i(\alpha)$ and $$H^{i(\alpha)}(F, S^{\alpha_1}\cQ_1\otimes \cdots \otimes S^{\alpha_{k+1}}\cQ_{k+1})\simeq S^{\overline{\alpha}}V.$$
\end{itemize}
\end{theo}

\subsection{A vanishing theorem}

As an application of the above theorem, we deduce a special vanishing theorem on flag varieties in this subsection, which will be one of the main ingredients of the proof of the main result. We also believe that this result has its own interests and may be useful in other contexts.

Let $\alpha_1, \ldots, \alpha_{k+1}$ be partitions of length $r_1,\ldots, r_{k+1}$ and let  $L=L_1^{\otimes a_1}\otimes\cdots\otimes L_k^{\otimes a_k}$ be an ample line bundle on $F$ such that $a_{i}-a_{i+1}\geq l>0$ for  $1\leq i\leq k-1$  and $a_k\geq l$ (or equivalently $L=L'^{\otimes l}\otimes L''$ where $L'$ and $L''$ are respectively  ample and nef line bundles  on $F$). For convention, we will write $a_{k+1}=0$.
Our goal is to compute the cohomology of vector bundles of the form
$$S^{\alpha_1}\cQ_1\otimes \cdots\otimes S^{\alpha_{k+1}}\cQ_{k+1}\otimes L.$$

 Let  $\tilde{\alpha}_i=\alpha_i+a_iI_i=( a_{i1}+a_i,\ldots,  a_{ir_i}+a_i)$ for $1\leq i\leq k+1$  and let $$\tilde{\alpha}=(\tilde{\alpha}_1, \tilde{\alpha}_2, \ldots, \tilde{\alpha}_{k+1})\in \mathbb Z^n$$ be a vector with $n$ entries. We then denote $$\tilde{\alpha}'=\tilde{\alpha}-(1,2, \ldots, n)=(\tilde{\alpha}_1', \ldots, \tilde{\alpha}_{k+1}'),$$ where $$\tilde{\alpha}_i'=\tilde{\alpha}_i-(r_1+r_2+\ldots+r_{i-1}+1, r_1+r_2+\ldots+r_{i-1}+2,\ldots, r_1+r_2+\ldots+r_{i-1}+r_i)$$ is the $i$-th block of $\tilde{\alpha}'$.  We write $\tilde{\alpha}_i'=(a'_{i1}, \ldots, a'_{ir_i})$, where $a'_{ij}=a_{ij}+a_i-\sum_{t=1}^{i-1}r_t-j$.

Note that $$S^{\alpha_1}\cQ_1\otimes \cdots\otimes S^{\alpha_{k+1}}\cQ_{k+1}\otimes L=S^{\tilde{\alpha}_1}\cQ_1\otimes \cdots\otimes S^{\tilde{\alpha}_{k+1}}\cQ_{k+1}.$$ Hence we need to study the number of inversions of $\tilde{\alpha}'$ to compute the cohomology of $$S^{\alpha_1}\cQ_1\otimes \cdots\otimes S^{\alpha_{k+1}}\cQ_{k+1}\otimes L.$$ By Theorem \ref{BBW}, we understand that if $\tilde{\alpha}'$ contains two identical entries, then each cohomology of the vector bundle in question is zero. Otherwise, we have the following vanishing theorem.
\begin{prop}\label{crucial}
In the above setting,
assume that any two entries of  $\tilde{\alpha}'$ are distinct,
there exist natural numbers $0\leq s_i\leq r_i$ for $2\leq i\leq k+1$ such that $$i(\tilde{\alpha})\leq \sum_{2\leq p\leq k+1}\sum_{1\leq i\leq s_j}a_{pi}-l(s_2+\ldots+s_{k+1})-\sum_{2\leq j\leq k+1}s_j^2.$$
 Thus, $H^q(F, S^{\alpha_1}\cQ_1\otimes \cdots\otimes S^{\alpha_{k+1}}\cQ_{k+1}\otimes L)=0$ for $$q>\sum_{2\leq p\leq k+1}\sum_{1\leq i\leq s_j}a_{pi}-l(s_2+\ldots+s_{k+1})-\sum_{2\leq j\leq k+1}s_j^2.$$

\end{prop}
\begin{proof}
Recall that $i(\tilde{\alpha})$ is the number of inversions of $\tilde{\alpha}'$.
We start to re-order  $\tilde{\alpha}'$. We  do it bloc by bloc. Note that each   $\tilde{\alpha}_i' $ is a strictly decreasing partition.

We first  assume that the first $s_2$ entries of $\tilde{\alpha}_2'$ need to move forward.
We may assume that the first $t_1$ entries of $\tilde{\alpha}_2'$ move forward $m_1$-steps to pass $\tilde{\alpha}_{1 t_1'}'$ and the next $t_2$  entries of $\tilde{\alpha}_2'$ move forward $m_2$-steps to pass $\tilde{\alpha}_{1 t_2'}'$, $\ldots$, and the  next $t_j$  entries of $\tilde{\alpha}_2'$ move forward $m_j$-steps to pass $\tilde{\alpha}_{1 t_j'}'$, where
$s_2=t_1+t_2+\cdots+t_j$ and $t_1'<t_2'<\cdots<t_j'$. Note that $m_1=r_1+1-t_1', \ldots, m_j=r_1+1-t_j'$.

Then we have the following inequalities:
\begin{eqnarray*}
&&a_{21}+a_2\geq a_{1t_1'}+m_1+t_1+a_1\\&&\;\;\;\vdots\\ && a_{2t_1}+a_2\geq a_{1t_1'}+m_1+t_1+a_1,\\
&&a_{2 (t_1+1)}+a_2\geq a_{1 t_2'}+t_1+m_2+t_2+a_1\\&&\;\;\;\vdots\\&& a_{2 (t_1+t_2)}+a_2\geq a_{1 t_2'}+t_1+m_2+t_2+a_1,\\
&&\vdots\\
&&a_{2 (t_1+\cdots+t_{j-1}+1)}+a_2\geq a_{1 t_j'}+(t_1+\cdots+t_{j-1}+m_j)+t_j+a_1\\&&\;\;\;\vdots\\&& a_{2 (t_1+\cdots+t_{j-1}+t_j)}+a_2\geq a_{1 t_j'}+(t_1+\cdots+t_{j-1}+m_j)+t_j+a_1.
\end{eqnarray*}
Since $a_1-a_2\geq l$, we take the sum of the above inequalities to get
$$\sum_{1\leq i\leq s_2}a_{2 i}\geq ls_2+\sum_{1\leq i\leq j}(t_im_i)+\sum_{1\leq i\leq j}t_i^2+\sum_{1\leq i_1<i_2\leq j}t_{i_1}t_{i_2}+\sum_{1\leq i\leq j} t_ia_{1t_i'}.$$
Note that $M_2:=\sum_{1\leq i\leq j}(t_im_i)$ is the number of inversions of this first step re-ordering. Moreover, since there are spaces for $t_i$ entries between $\tilde{a}_{1(t_i'-1)}' $ and $\tilde{a}_{1 t_i'}'$, we have $$a_{1 t_{i-1}'}\geq a_{1(t_i'-1)}\geq a_{1 t_i'}+t_i$$ for $2\leq i\leq j$. Thus $$\sum_{1\leq i\leq j} t_ia_{1t_i'}\geq s_2a_{1t_j'}+\sum_{1\leq i_1<i_2\leq j}t_{i_1}t_{i_2}.$$ We finally have $$\sum_{1\leq i\leq s_2}a_{2 i}\geq ls_2+ s_2a_{1t_j'}+s_2^2+M_2\geq ls_2+s_2^2+M_2.$$

We now assume by induction that after the re-ordering of the first $h$ block, there exist natural numbers $0\leq s_i\leq r_i$ for $2\leq i\leq h$ such that
$$\sum_{2\leq p\leq h}\sum_{1\leq i\leq s_2}a_{p i}\geq l(s_2+\cdots+s_h)+s_2^2+\cdots+s_h^2+M_2+\cdots+M_h,$$ where $M_i$ is the number of inversions while moving the entries of the  $i$-th bloc forward.

After the re-ordering of the first $h$ blocs we have a partition of length $r_1+r_2+\cdots+r_h$. We denote this partition with strictly decreasing numbers by $\mu=(c_1, c_2 \ldots, c_{r_1+\cdots+r_h})$.  We need to compare this partition with the original partition. We thus  write $$(\alpha_1, \ldots, \alpha_h)=(b_1, b_2, \ldots, b_{r_1+\cdots+r_h})\in \mathbb Z^{r_1+r_2+\cdots+r_h}.$$  Let $f: \{1, 2, \ldots,n\}\rightarrow \{1,\ldots, k\}$ be the map such that the $i$-th entry of $$(a_1I_1, a_2I_2,\ldots, a_kI_k, a_{k+1}I_{k+1})$$ is $a_{f(i)}$.

After reviewing the re-ordering we have done, one sees that either $c_j=b_{j'}-j'+a_{f(j')}$ for some $j'\leq j$ or $c_j=b_{j''}-j''+a_{f(j'')}$ for some $j''>j$ and in the latter case $b_{j''}-j''+a_{f(j'')}>b_{j'}-j'+a_{f(j')}$ for some $j'\leq j$.

We can now start the re-ordering for the $(h+1)$-th bloc.  For the  vector
$$(\mu, \tilde{\alpha}_{h+1}')\in \mathbb Z^{r_1+r_2+\cdots+r_{h+1}},$$ we  assume that the first $t_1$ entries of $\tilde{\alpha}_{h+1}'$ move forward $m_1$-steps to pass $c_{q_1}$ and the next $t_2$  entries of $\tilde{\alpha}_{h+1}'$ move forward $m_2$-steps to pass $c_{q_2}$, $\ldots$, and the  next $t_j$  entries of $\tilde{\alpha}_{h+1}'$ move forward $m_j$-steps to pass $c_{q_j}$, where
$s_{h+1}=t_1+t_2+\cdots+t_j$. Similarly, we have in this case $r_1+r_2+\cdots+r_h+1-q_1=m_1, \ldots, r_1+r_2+\cdots+r_h+1-q_j=m_j$. Then
\begin{eqnarray*}
a_{(h+1)1}+a_{h+1}\geq \cdots \geq a_{(h+1)t_1}+a_{h+1}&&\geq (c_{q_1}+q_1)+m_1+t_1,\\
a_{(h+1) (t_1+1)}+a_{h+1}\geq\cdots\geq a_{(h+1) (t_1+t_2)}+a_{h+1}&&\geq (c_{ q_2}+q_2)+(t_1+m_2)+t_2,\\
&&\vdots\\
a_{(h+1)(t_1+\cdots+t_{j-1}+1)}+a_{h+1}\geq \cdots\geq a_{(h+1)(t_1+\cdots+t_{j-1}+t_j)}+a_{h+1}&&\geq (c_{ q_j}+q_j)+\\&&(t_1+\cdots+t_{j-1}+m_j)+t_j.
\end{eqnarray*}

Similarly, we have $c_{q_{i-1}}+q_{i-1}\geq (c_{q_i}+q_i)+t_i$. Summing up, we have
$$\sum_{1\leq i\leq s_{h+1}}a_{(h+1)i}\geq s_{h+1}(c_{q_j}+q_j-a_{h+1})+M_{h+1}+s_{h+1}^2.$$

 By the above discussion about $c_j$, we always have $c_{q_j}\geq b_{q_j'}-q_j'+a_{f(q_j')} $ for some $q_j'\leq q_j$. We also have $f(q_j')\leq a_h$. Thus $c_{q_j}+q_j\geq a_{f(q_j')}\geq a_{h+1}+l$. Hence
 $$\sum_{1\leq i\leq s_{h+1}}a_{(h+1)i}\geq ls_{m+1}+M_{h+1}+s_{h+1}^2.$$

 We then finish the induction and conclude that
 $$\sum_{2\leq p\leq k+1}\sum_{1\leq i\leq s_p}a_{p i}\geq l(s_2+\cdots+s_{k+1})+s_2^2+\cdots+s_{k+1}^2+M_2+\cdots+M_{k+1}.$$
  Since $i(\tilde{\alpha})$, which is the number of inversion of $\tilde{\alpha}'$, is exactly $M_2+\cdots+M_{k+1}$. Therefore, combining Theorem \ref{BBW}, we conclude the proof.

\end{proof}
\subsection{A Filtration on the kernel bundle}

As we have explained in the introduction, the property of the kernel bundles plays an important role in the study of syzygies. We recall in this subsection the work of Manivel in \cite{Man} about the kernel bundles on flag varieties. Manivel applied Schur functors to give  natural filtrations on the kernel bundles so that the associated graded vector bundles admit projective resolutions consisting of tensor products of Schur powers of the tautological sub and quotient bundles.

In this subsection, we still work on $ F:=\mathrm{Fl}(n_k, \ldots, n_1;V)$.
Let $a_1\geq a_2\geq\cdots \geq a_k\geq 0$ be a sequence of integer and let $L=L_1^{\otimes a_1}\otimes\cdots \otimes L_k^{\otimes a_k}$. We denote by $$M_L:=\mathrm{Ker}(H^0(F, L)\otimes \cO_F\rightarrow L)$$ the kernel bundle of $L$.

 We fix some notations. Let $\alpha_i=a_iI_i=(a_i,\ldots,a_i)$ be a partition of length $r_i$.  Let $$\alpha=(\alpha_1,\ldots,\alpha_k, \overbrace{0,\ldots,0}^{r_{k+1}})$$ be a partition of length $n$.
We now define $\alpha_{\geq i}$ a partition of length $r_1+\cdots +r_{i+1}$ for $1\leq i\leq k$: $$\alpha_{\geq i}:=(\alpha_1,\ldots, \alpha_i,\overbrace{0,\ldots,0}^{r_{i+1}})-(a_{i+1}I_{1},\ldots, a_{i+1}I_{i},\overbrace{0,\ldots,0}^{r_{i+1}}).$$ If $a_i>a_{i+1}$,  there are exactly $r_{i+1}$ zeros among the  entries  of $\alpha_{\geq i}$. We also define $\tilde{L}_i=(L_1\otimes\cdots \otimes L_i)^{\otimes a_i}\otimes L_{i+1}^{\otimes a_{i+1}}\otimes\cdots\otimes L_k^{\otimes a_k}$.

 Since $L=S^{\alpha_1}\cQ_1\otimes \cdots\otimes S^{\alpha_k}\cQ_k$, $H^0(F, L)$ can be naturally identified with $S^{\alpha}V$ and the evaluation map is identified with $\mathrm{ev}_{\alpha}: S^{\alpha}\Sigma_0\rightarrow L.$

The map $\mathrm{ev}_{\alpha}$ factors naturally through the following sequence of    Schur powers of quotient maps
 \begin{eqnarray*}\mathrm{ev_k}: S^{\alpha}\Sigma_0 &&\rightarrow S^{\alpha}(\Sigma_0/\Sigma_k)=S^{\alpha_{\geq k-1}}(\Sigma_0/\Sigma_k)\otimes \tilde{L}_k,\\
& \vdots &\\
\mathrm{ev}_i:S^{\alpha_{\geq i}}(\Sigma_0/\Sigma_{i+1})\otimes\tilde{L}_{i+1}&&\rightarrow S^{\alpha_{\geq i}}(\Sigma_0/\Sigma_{i})\otimes\tilde{L}_{i+1}=S^{\alpha_{\geq i-1}}(\Sigma_0/\Sigma_{i})\otimes\tilde{L}_{i},\\
&\vdots&\\
\mathrm{ev}_1: S^{\alpha_{\geq 1}}(\Sigma_0/\Sigma_{2})\otimes\tilde{L}_{2} &&\rightarrow S^{\alpha_{\geq 1}}(\Sigma_0/\Sigma_{1})\otimes\tilde{L}_{2}=\tilde{L}_1=L.
 \end{eqnarray*}

 Therefore, $\mathrm{ker}(\mathrm{ev_k})\subset \mathrm{ker}(\mathrm{ev_{k-1}}\circ \mathrm{ev_k})\subset\cdots \subset  \mathrm{ker}(\mathrm{ev_{1}}\circ\cdots\circ \mathrm{ev_k})=M_L$ gives a filtration of $M_L$, with the graded bundles $N_i:=\mathrm{ker} (\mathrm{ev}_i)$ for $1\leq i\leq k$.

The following lemma  (see \cite[Proposition in Subsection 2.2]{Man}) is crucial.

 \begin{lem}\label{weyman}
There are exact Schur complexes:
 \begin{eqnarray*}
\cdots \rightarrow C(i)_j\rightarrow \cdots \rightarrow C(i)_1\rightarrow N_i\rightarrow 0,
 \end{eqnarray*}
 where $$C(i)_j=\big[\bigoplus_{|\rho|=j, |\nu|=|\alpha_{\geq i}|-j}  \big(S^{\rho}\cQ_{i+1}\otimes S^{\nu}(\Sigma_0/\Sigma_{i+1})\big)^{\oplus c_{\rho*, \nu}^{\alpha_{\geq i}}}\big]\otimes \tilde{L}_{i+1}$$ and $c_{\rho^*, \nu}^{\alpha_{\geq i}}$ is a Littlewood-Richardson coefficient: which is the multiplicity of $S^{\alpha_{\geq i}} W$ in the tensor product $S^{\rho^*}W\otimes S^{\nu} W$, where $\rho^*$ is the conjugate of $\rho$ and $W$ is a vector space of dimension large enough.
 \end{lem}

 We also need the following vanishing theorem, which is a  generalization of \cite[Corollary in Section 3]{Man}, to prove the main theorem.

 Let $\rho=\{\rho_1,\ldots, \rho_k\}$ be a set of $k$ partitions and define $S^{\rho}\cF:=S^{\rho_1}(\Sigma_0/\Sigma_1)\otimes\cdots\otimes S^{\rho_k}(\Sigma_0/\Sigma_k)$. Recall that $\Sigma_0=V\otimes \cO_F$ is the trivial vector bundle. Thus the quotient bundles $\Sigma_0/\Sigma_i$ are nef. We observe that the line bundles $\tilde{L}_i$ in Lemma \ref{weyman} can also be expressed as a Schur power $S^{\rho}\cF$.

 \begin{prop}\label{finalvanishing} Given partitions $\beta_i=(b_{i1}, \cdots, b_{ir_i})$ for $1\leq i\leq k+1$ and let $L$ be an ample line bundle on $F$ which can be expressed as $L'^{\otimes l}\otimes L''$ with $L'$ ample and $L''$ nef. Assume that $$q>\sum_{2\leq i\leq k+1}\sum_{1\leq j\leq s_i}b_{ij}-l(\sum_{2\leq i\leq k+1}s_i)-\sum_{2\leq i\leq k+1}s_i^2$$ for any integers $0\leq s_2\leq r_2, \ldots, 0\leq s_{k+1}\leq r_{k+1}$.
 Then \begin{eqnarray}\label{finalvanishing1}H^q(F, S^{\rho}\cF\otimes  S^{\beta_1}\cQ_1\otimes \cdots\otimes S^{\beta_{k+1}}\cQ_{k+1} \otimes L)=0,\end{eqnarray} for any partitions $\rho_1,\ldots, \rho_k$.
 \end{prop}
\begin{proof}
  We combine the same argument in the proof of  \cite[Corollary in Subsection 3.3]{Man} with Proposition \ref{crucial} to give a proof.  We run inductions on $|\rho|:=\sum_{1\leq i\leq k}|\rho_i|$. When $|\rho|=0$, this is just Proposition \ref{crucial}.  We now assume that the vanishing (\ref{finalvanishing1}) holds when $|\rho|<M$. We now assume that $|\rho|=\sum_{1\leq i\leq k}|\rho_i|=M$. Then there exists a set $\rho'$ of $k$ partitions such that $|\rho'|=M-1$ and $S^{\rho}\cF$ is a direct summand of $S^{\rho'}\cF\otimes (\Sigma_0/\Sigma_i)$ for some $1\leq i\leq k$. We have the short exact sequence:
\begin{eqnarray*}0\rightarrow S^{\rho'}\cF\otimes \Sigma_i\rightarrow  S^{\rho'}\cF\otimes \Sigma_0\rightarrow  S^{\rho'}\cF\otimes (\Sigma_0/\Sigma_i)\rightarrow 0.
\end{eqnarray*}
Then by induction, for $$q>\sum_{2\leq i\leq k+1}\sum_{1\leq j\leq s_i}b_{ij}-l(\sum_{2\leq i\leq k+1}s_i)-\sum_{2\leq i\leq k+1}s_i^2$$ for any $0\leq s_2\leq r_2, \ldots, 0\leq s_{k+1}\leq r_{k+1}$, we have
$$H^q(F, S^{\rho'}\cF\otimes \Sigma_0\otimes  S^{\beta_1}\cQ_1\otimes \cdots\otimes S^{\beta_{k+1}}\cQ_{k+1} \otimes L)=0,$$because $\Sigma_0=V\otimes \cO_F$ is the trivial vector bundle. We also have
 $$H^{q+1}(F, S^{\rho}\cF\otimes \Sigma_i\otimes  S^{\beta_1}\cQ_1\otimes \cdots\otimes S^{\beta_{k+1}}\cQ_{k+1} \otimes L)=0,$$ since there exists a filtration of $\Sigma_i$  whose graded vector bundles are $\cQ_{i+1}, \ldots, \cQ_{k+1}$ and the induction implies that $$H^{q+1}(F, S^{\rho'}\cF\otimes \Sigma_i\otimes  S^{\beta_1}\cQ_1\otimes \cdots \otimes S^{\beta_j'}\cQ_j\otimes\cdots \otimes S^{\beta_{k+1}}\cQ_{k+1} \otimes L)=0,$$ where $S^{\beta_j'}\cQ_j$ is any direct summand of $S^{\beta_j}\cQ_j\otimes \cQ_j$ for $i+1\leq j\leq k+1$.

  Hence $$H^q(F, S^{\rho}\cF\otimes  S^{\beta_1}\cQ_1\otimes \cdots\otimes S^{\beta_{k+1}}\cQ_{k+1} \otimes L)=0.$$
\end{proof}
\section{Rational homogeneous varieties}
We recall some basic facts about rational homogeneous varieties in this section (see \cite{Man1, LM} and  references therein).

A classical result of Borel and Remmert says that a homogeneous variety is  a direct product of an abelian variety with a rational homogeneous variety. A rational homogeneous variety can be written as $G/P$, where $G$ is a semisimple algebraic group and $P$ is a parabolic subgroup. Furthermore, $G/P$ can be decomposed as $G_1/P_1\times \cdots \times G_k/P_k$, where $G_i$ are simple algebraic groups.

For a simple algebraic group $G$, the homogeneous spaces $G/P$ can be characterized by its Lie algebra $\mathfrak{g}$ and the Dynkin diagram $\Delta$. The Dynkin diagrams of classical groups correspond to type $A_n$, $B_n$, $C_n$, and $D_n$. The exceptional ones are type $E_6$, $E_7$, $E_8$, $F_4$, and $G_2$.
The isomorphic classes of $G$-homogeneous variety are in bijective with the finite subsets of nodes of $\Delta$.

Fix a Borel subgroup $B$ of $G$ and a maximal torus. Let $\Delta=\{\alpha_1, \cdots, \alpha_n\}$ be the set of simple roots and the fundamental weights are $\lambda_1,\cdots, \lambda_n$. The standard parabolic subgroups containing $B$ is determined by a subset $S$ of $\Delta$. For each $S$,  the corresponding parabolic subgroup $P_S$ is characterized by $\alpha_i\notin \mathrm{Lie}(P)$ for $i\in S$. By a classical result of Chevalley, we  also have $$\Pic(G/P_S)\simeq P(S):=\{\sum_{\alpha_i\in S}n_i\lambda_i\mid n_i\in \mathbb Z\},$$ where $P(S)$ is the sub-lattice generated by the weights $\lambda_i$ dual to the roots $\alpha_i\in S$. Indeed, for each weight $\lambda\in P(S)$, there exists a corresponding one-dimensional representation $L_{\lambda}$ of $P_S$ and we have a line bundle $\cL_{\lambda}:=G\times^{P_S}L_{\lambda}$ over $G/P_S$.

As explained in the introduction, we focus on rational homogeneous varieties of type $B_n$, $C_n$, $D_n$, and $G_2$. One of the main reason is that they are naturally closed subvarieties of flag varieties and  their Picard groups can also be described with the help of the Picard goups of flag varieties.

  \begin{lem}
 Assume that $X$ is a homogeneous variety of type $C_n$, or $B_n$, or $D_n$, or $G_2$, then $X$ is a closed subvariety of a flag variety $F=\mathrm{Fl}(n_k,\ldots,n_1; V)$ and can be identified with a connected component of the zero locus of a global section of certain homogeneous vector bundle $W^*$ and $\codim_FX=\rank W^*$. In the $C_n$ case, $W=\wedge^2\Sigma_1$, in the $B_n$ or $D_n$ case, $W=S^2\Sigma_1$. In the $G_2$ case, $W$ will be explained in Subsection 3.4.

 \end{lem}
This lemma is certainly well known. We simply recall explicitly the structures of rational homogeneous varieties of type $C_n$, or $B_n$, or $D_n$, or $G_2$ in the following subsections.

\subsection{Rational homogeneous varieties of type $C_n$}

Let $V$ be a linear space of dimension $2n$ with $n\geq 3$ and $\omega\in\wedge^2V^*$ be a non-degenerate $2$-from. The  symplectic group $\mathrm{Sp}(n)=\mathrm{Sp}(V, \omega)$ is of type $C_n$. Thus the rational homogeneous varieties of type $C_n$ are isotropic flag varieties. Fixing  a symplectic vector space $(V, \omega)$  of dimension $2n$, for any sequence $1\leq n_k<\cdots<n_1\leq n$, let $\mathrm{SFl}(n_k,\ldots,n_1;V)$ be the flag variety parametrizing  subspaces $V_k\subset \cdots\subset V_1$ such that $\dim V_i=n_i$ and $\omega|_{V_1}=0$.   It is clear that $\mathrm{SFl}(n_k,\ldots,n_1; V)$ is a closed subvariey of   $\mathrm{Fl}(n_k,\ldots, n_1;V)$ defined by the global section $\wedge^2\Sigma_1^*$ corresponding to $\omega$. By comparing the root lattices, the restriction map $$\Pic(\mathrm{Fl}(n_k,\ldots, n_1; V))\rightarrow \Pic(\mathrm{SFl}(n_k,\ldots, n_1; V))$$ is an isomorphism.

\subsection{Rational homogeneous varieties of type $B_n$}
Let $V$ be a linear space with a non-degenerate quadratic form $q\in S^2V^*$. The group $\mathrm{SO}(V, q)$ is of type  $B_n$ (resp. $D_n$) when $\dim V=2n+1$ with $n\geq 2$ (resp. $\dim V=2n$ with $n\geq 4$). We have a similar description of the homogeneous varieties of type $B_n$ or $D_n$ as the orthogonal flag varieties.  But we need to be a bit more careful.

In the $B_n$ case, the orthogonal Grassmannian $\mathrm{OGr}(n;V)$ parametrizing maximal isotropic subspaces  is connected and the Pl\"ucker divisor of $\mathrm{Gr}(n, V)$ restricted to $\mathrm{OGr}(n; V) $ is twice the generator of the Picard group of the latter. A rational homogeneous variety of
type $B_n$ is  an isotropic flag varieties $\mathrm{OFl}(n_k,\ldots,n_1; V)$
parametrizing flags $V_k\subset\cdots\subset V_1$ with $\dim V_i=n_i$ and $V_1$ is isotropic with respect to $q$. Note that $\mathrm{OFl}(n_k,\ldots,n_1; V)$ is defined as the zero loci of the global section of $S^2\Sigma_1^*$   corresponding to $q$. The injective restriction map $$\Pic(\mathrm{Fl}(n_k,\ldots,n_1; V))\rightarrow \Pic(\mathrm{OFl}(n_k,\ldots,n_1; V))$$ is an isomorphism when $n_1<n$ and the image is a sub-lattice of  index $2$ when $n_1=n$.

\subsection{Rational homogeneous varieties of type $D_n$}

Fix a linear space $V$ of dimension $2n$ with a non-degenerate quadric form $q$, the orthogonal Grassmannian $\mathrm{OGr}(n; V)$ parametrizing maximal isotropic subspaces has two (isomorphic) connected component $\mathrm{OGr}(n; V)^{+}$ and $\mathrm{OGr}(n; V)^{-}$.    Moreover, the Pl\"ucker divisor of $\mathrm{Gr}(n, V)$ restricted to one of the component of  $\mathrm{OGr}(n;V) $ is again twice the generator of the Picard group of the latter.
 Indeed, taking a general hyperplane $V'$ of $V$, it is easy to see that any connected component of $\mathrm{OGr}(n;V)$ is isomorphic to $\mathrm{OGr}(n-1;V')$.
Given a isotropic subspace $W$ of dimension $n-1$ of $V$, there are exactly two maximal isotropic subspaces $W_1$ and $W_2$ containing $W$, which belong to different components of $\mathrm{OGr}(n; V)$. A rational homogeneous variety $X$ of type $D_n$ is  one of the following: for $0<n_k<\cdots<n_1\leq n-2$, $\mathrm{OFl}(n_k,\ldots,n_1; V)$, $\mathrm{OFl}(n_k,\ldots,n_1,n; V)^{\epsilon}$ with $\epsilon=+$ or $-$,  and  $$\mathrm{OFl}(n_k,\ldots, n_1,n-1; V)=\mathrm{OFl}(n_k, \ldots, n-1,\ldots,n; V)^{+}\times_{\mathrm{OFl}(n_k,\ldots,n_1; V)}\mathrm{OFl}(n_k, \ldots, n-1,\ldots,n; V)^{-}.$$

 For $0<n_k<\cdots<n_1\leq n-2$ and the natural embedding $$ \mathrm{Fl}(n_,k\ldots, n_1; V)\hookrightarrow \mathrm{Fl}(n_k,\ldots, n_1; V),$$ the restriction map
$\Pic(\mathrm{Fl}(n_k,\ldots,n_1; V))\rightarrow \Pic(\mathrm{OFl}(n_k,\ldots,n_1; V))$ is  an isomorphism. For the natural embedding $$ \mathrm{OFl}(n_k,\ldots, n_1,n; V)^{\epsilon}\hookrightarrow \mathrm{Fl}(n_k,\ldots, n_1,n; V),$$ the image of the injective restriction map
$\Pic(\mathrm{Fl}(n_k,\ldots, n_1,n; V))\rightarrow \Pic(\mathrm{OFl}(n_k,\ldots, n_1,n; V)^{\epsilon})$ is a sub-lattice of index $2$.  For the natural embedding $$ \mathrm{OFl}(n_k,\ldots, n_1, n-1; V)\hookrightarrow \mathrm{Fl}(n_k,\ldots, n_1, n-1; V),$$ the restriction map $\Pic(\mathrm{Fl}(n_k,\ldots, n_1, n-1; V))\rightarrow \Pic(\mathrm{OFl}(n_k,\ldots, n_1, n-1; V))$ is a  embedding between lattices whose cokernel is a  lattice of rank $1$.

\subsection{Rational homogeneous varieties of type $G_2$}

 The rational homogeneous varieties of type $G_2$ are also well understood (see \cite[Subsectin 6.1]{LM} and \cite[Appendix B]{K1}). Fix a vector space $V$ of dimension $7$. The group $G_2$  can  be realized as the subgroup of $\mathrm{GL}(V)$ fixing  a generic element $\lambda$ of $  \wedge^3V^*$. The group $G_2$ fixes a non-degenerate symmetric two form $\theta$ on $V$ and hence $G_2$ is also a subgroup of $O(V, \theta)$. The three non-trivial homogeneous varieties of type $G_2$ are: the quadric $Q$ in $\mathbb P(V)$ corresponding to $\theta$, $ X:=G_2\mathrm{Gr}(2, V)$ which is the zero locus of the global section corresponding to  $\lambda$
 of $\mathcal{Q}^*(1)$  on $\mathrm{Gr}(2, V)$, where $\mathcal Q$ is the tautological quotient bundle on $\mathrm{Gr}(2, V)$ \footnote{Note that
 $\wedge^3V^*\simeq \wedge^4V=\Gamma(\mathrm{Gr}(2, V), \mathcal{Q}^*(1))$.},   and $\mathbb P_X(\Sigma)$, where $\Sigma$ is the tautological sub-bundle on $\mathrm{Gr}(2, V)$. We remark that $X$ is a Fano manifold of dimension $5$ and of index $3$. We have the following commutative diagram:
 \begin{eqnarray*}
 \xymatrix{
X \ar@{^{(}->}[d]& \mathbb P_X(\Sigma)\ar[l]\ar@{^{(}->}[d]\ar[r] & Q\ar@{=}[d]\\
\mathrm{OGr}(2, V)\ar@{^{(}->}[d]&\mathrm{OFl}(1, 2; V)\ar[l]\ar@{^{(}->}[d]\ar[r] & Q\ar@{^{(}->}[d]\\
\mathrm{Gr}(2, V)&\mathrm{Fl}(1, 2; V)\ar[l]^p\ar[r]& \mathbb P(V).}
 \end{eqnarray*}
 The varieties in the first line are exactly the homogenous varieties of type $G_2$. The commutative diagrams in the left are Cartesians. Thus $\mathbb P_X(\Sigma)$ is the zero locus of a global section of $p^*(\cQ(1))$ on $\mathrm{Fl}(1, 2; V)$. It is also easy to verify that the restriction maps on Picard groups via the vertical morphisms in the above diagram are always bijective.

\section{Plethysm}
\subsection{Some results on partitions}

The general plethysm is still an open problem. But the Schur decompositions of $\wedge^j(\wedge^2 V)$ or $\wedge^j(S^2V)$ are known (see \cite[Subsection 1.8]{Mac} or \cite[2.3.9]{W}). Given a partition $\alpha=(a_1,\ldots, a_n)$ of a dimension $n$ space $V$, we consider its Young diagram and let $r_{\alpha}$ be the length of the diagonal of its Young diagram, called the rank of $\alpha$. Let $\alpha^*=(a_1^*,\ldots, a^*_n)$ be the conjugate of $\alpha$, namely the rows of the Young diagram of $\alpha^*$ are the columns of the Young diagram of $\alpha$. We can also express $\alpha$ as $(\mu_1, \ldots, \mu_r|\rho_1,\ldots,\rho_r)$, where $r=r_{\alpha}$, $\mu_i=a_i-i$ and $\rho_i=a^*_i-i$.

Then \begin{eqnarray}\label{plethysm1}\wedge^j(\wedge^2 V)=\bigoplus_{\alpha}S^{\alpha}V,\end{eqnarray} where $\alpha$ goes through all partitions of weight $2j$, which can be expressed as $(\mu_1, \ldots,\mu_r|\mu_1+1,\ldots,\mu_r+1)$ where $r$ is the rank.

Similarly, we have \begin{eqnarray}\label{plethysm2}\wedge^j(S^2 V)=\bigoplus_{\alpha}S^{\alpha}V,\end{eqnarray} where $\alpha$ goes through the conjugate of the partitions in (\ref{plethysm1}).
\begin{lem}\label{precise}
Let $V$ be a $n$-dimensional vector space and $\alpha=(a_1, \ldots, a_n)$ be a partition.
\begin{itemize}
\item If $S^{\alpha}V$ is a direct summand of $\wedge^j(\wedge^2V)$, we have $$\sum_{1\leq i\leq s}a_i\leq j+\frac{s^2-s}{2},$$ for any $1\leq s\leq n$.
 \item If $S^{\alpha}V$ is a direct summand of  $\wedge^j(S^2V)$, we have $$\sum_{1\leq i\leq s}a_i\leq  j+\frac{s^2+s}{2},$$ for any $1\leq s\leq n$.
\end{itemize}
\end{lem}
\begin{proof}
Assume that  $S^{\alpha}V$ is a  direct summand of $\wedge^j(\wedge^2V)$. Then by (\ref{plethysm1}), $\alpha$ can be expressed as $(\mu_1, \ldots, \mu_r|\mu_1+1,\ldots, \mu_r+1)$, where $r$ is the rank of $\alpha$ and $\mu_i=a_i-i$ for $1\leq i\leq r$.  Note that $2j=|\alpha|\geq r^2$.

We claim that $\sum_{1\leq i\leq r}a_i=j+\frac{r^2-r}{2}$. Indeed,
\begin{eqnarray*}
2j=|\alpha |=\sum_{1\leq i\leq n}a_i= \sum_{1\leq i\leq r}a_i+\sum_{r+1\leq i\leq n}a_i
\end{eqnarray*}
From the expression of $\alpha$, we see that the Young diagram of $\alpha$ has the first $r$ columns with respectively $\mu_1+2, \mu_2+3, \ldots, \mu_r+r+1$ elements. Thus $$\sum_{r+1\leq i\leq n}a_i=\sum_{1\leq i\leq r}(\mu_i+i+1-r)=\sum_{1\leq i\leq r}a_i-r(r-1).$$ Thus $\sum_{1\leq i\leq r}a_i=j+\frac{r^2-r}{2}$.

If $s<r$, we   have \begin{eqnarray*}&&\sum_{1\leq i\leq s}a_i= \sum_{1\leq i\leq r}a_i-\sum_{s+1\leq i\leq r}a_i\\&&\leq j+\frac{r^2-r}{2}-\sum_{s
+1\leq i\leq r}r=j+\frac{r^2-r}{2}-r^2+rs\leq j+\frac{s^2-s}{2}.
\end{eqnarray*}
If $s>r$, we have \begin{eqnarray*}&&\sum_{1\leq i\leq s}a_i= \sum_{1\leq i\leq r}a_i+\sum_{r+1\leq i\leq s}a_i\\&&\leq j+\frac{r^2-r}{2}+\sum_{r
+1\leq i\leq s}r\leq j+\frac{s^2-s}{2} .
\end{eqnarray*}

Assume that $S^{\alpha}V$ is a direct summand of $\wedge^j(S^2V)$. Then for the conjugate $\alpha^*$, $S^{\alpha^*}V$ appears as a direct summand of $\wedge^i(\wedge^2V)$ by (\ref{plethysm1}) and  (\ref{plethysm2}). Let $r$ be the rank of $\alpha$, then $\sum_{1\leq i\leq r}a_i=j+\frac{r^2-r}{2}+r=j+\frac{r^2+r}{2}$. Then, we conclude similarly that if $s<r$, $$\sum_{1\leq i\leq s}a_i \leq j+\frac{s(s+1)}{2},$$ and if $s>r$, $$\sum_{1\leq i\leq s}a_i \leq j+\frac{s(s+1)}{2}.$$
 \end{proof}

\subsection{Weights of the graded bundles}

We have explained that rational homogeneous varieties of type $B_n$, $C_n$, $D_n$, and $G_2$ can be realized as  a connected component of the zero locus of a global section of a homogeneous bundle $W^*$ on a flag variety $F$.

We write $$\wedge^jW=\bigoplus_{\alpha} S^{\alpha}\Sigma_1,$$ where in the type $\text{C}$ case, $\alpha$ goes through all  partitions of weight $2j$ of length $\leq n_1$, which can be expressed as $(\mu_1, \ldots, \mu_r|\mu_1+1,\ldots, \mu_r+1)$ and in the type $\text{B}$ or $\text{D}$ case,  $\alpha$ goes through all  partitions of weight $2j$, of length $\leq n_1$, which can be expressed as $(\mu_1+1, \ldots, \mu_r+1|\mu_1, \ldots, \mu_r)$.

Since we have the filtration $\Sigma_1\supset \Sigma_2\supset \cdots\supset \Sigma_k\supset 0$ with successive graded bundles $\cQ_2,\ldots, \cQ_{k+1}$, $S^{\alpha}\Sigma_1$ also admits a filtration whose graded bundles are \begin{eqnarray}\label{quotients}S^{\rho_2}\cQ_2\otimes \cdots\otimes S^{\rho_{k+1}}\cQ_{k+1},\end{eqnarray} where $\rho_i$, $2\leq i\leq k+1$ go through all partitions such that $\sum_{2\leq i\leq k+1}|\rho_i|=|\alpha|=2j$ and $S^{\rho_2}D\otimes \cdots\otimes S^{\rho_{k+1}}D$ contains $S^{\alpha}D$ as a direct summand for some vector space $D$ of sufficiently large dimension (see for instance \cite[Exercise 6.11]{FH}).

\begin{lem}\label{precise1}Let $\rho=(\rho_2, \ldots, \rho_{k+1})$ be as in (\ref{quotients}). For any $s$ entries $a_1,\ldots, a_s$ of $\rho$,
$$\sum_{1\leq i\leq s}a_i\leq j+\frac{s(s-1)}{2}$$ in the type $C$ cases and
$$\sum_{1\leq i\leq s}a_i\leq j+\frac{s(s+1)}{2}$$ in the type $B$ or $D$ cases.
\end{lem}
\begin{proof}We know by   applying repeatedly  the Littlewood-Richardson rule or the Pieri rule that $\sum_{1\leq i\leq s}a_i$ is smaller or equal to the sum of the first $s$ entries of $\alpha$ and we then conclude by Lemma \ref{precise}.
\end{proof}
\section{Syzygies of Isotropic flag varieties}

Let $V$ be a vector space  and $F:=\mathrm{Fl}(n_k,\ldots, n_1; V)$ be a flag variety. We have seen that a  rational homogeneous variety $X$ of type $B_n$, $C_n$, or $D_n$ can be realized as a connected component of the zero locus $\tilde{X}$ of a global section of a homogeneous vector bundle $W^*$ on $F$ and $W=\wedge^2\Sigma_1$ in the type $C_n$ case and $W=S^2\Sigma_1$ in the   $B_n$ or $D_n$ type case.

 We then have a long exact sequence:
\begin{eqnarray}\label{resolution}0\rightarrow \wedge^{\rank W}W\rightarrow \cdots \wedge^2W\rightarrow\cO_{F}\rightarrow \cO_{\tilde X}\rightarrow 0.
\end{eqnarray}

 We will also denote by $q: F=\mathrm{Fl}(n_k,\ldots, n_1; V)\rightarrow G:=\mathrm{Gr}(n_1, V)$ the natural morphism from the flag variety to the Grassmannian. Then $\Sigma_1$ on $F$ is also the pull-back of the tautological subbundle $\Sigma_G$ on $G$ and $q$ is indeed the $\mathrm{Fl}(n_k,\ldots, n_2; \Sigma_G)$-bundle over $G$. We will write $L_G$ the Pl\"ucker divisor on $G$.

We write $W_G=\wedge^2\Sigma_G$ or $S^2\Sigma_G$ respectively in the type $C_n$ case or the type $B_n$ or $D_n$ case. We denote by $\alpha$ the partition $(a_1I_1,\ldots,  a_{k},\ldots, a_kI_k,0\ldots,0)$ and $\tilde{\alpha}=(a_2I_2,\ldots, a_{k}I_k, 0\ldots,0)$, where $I_i$ are as in Section 2.
\begin{lem}\label{restriction}Under the above setting, the restriction map $H^0(F, L)\rightarrow H^0(X, L|_X)$ is surjective for any ample line bundle $L$ on $F$.
\end{lem}

\begin{proof}

We know that $H^0(X, L)$ is an irreducible representation under the action of the symplectic or the orthogonal group and the restriction map is equivariant under the group action. Thus it suffices to show that the restriction map is nonzero. This follows easily from the fact that $L$ is globally generated.
\end{proof}
\begin{rema}
Since any ample line bundle $L$ on $F$ satisfies property $(N_0)$, this implies that $L|_{X}$ also satisfies Property $(N_0)$, which  is indeed  a classical result.
\end{rema}
\subsection{Symplectic flag varieties}

\begin{theo}\label{symplectic}

 Assume that $X=\mathrm{SFl}(n_k,\ldots, n_1; V)$ is a symplectic flag variety.
Let $L$ be a line bundle on $X$, which can be expressed as $H^{\otimes l}\otimes M$, where $H$ is an ample line bundle and $M$ is a nef line bundle. Then
\begin{itemize}
\item[(1)] When $k\leq 2$, then $L$ satisfies Property $(N_p)$ when $l\geq p$;
\item[(2)] When $l\geq p\geq \frac{n_1}{2}-1$, $L$ satisfies Property $(N_p)$;
\item[(3)] In general, when $l\geq \mathrm{max}\{p, \frac{p+1}{n_1}+\frac{n_1-3}{2}\}$, $L$ satisfies Property $(N_p)$.
\end{itemize}

 \end{theo}

\begin{proof}
Consider the natural embedding $\iota: X\hookrightarrow F:=\mathrm{Fl}(n_k,\ldots, n_1; V)$. Since the restriction map induces an isomorphism between the Picard groups, $L$ is the restriction on $X$ of a line bundle  on $F$.

For simplicity of the notations, we now write $L=L_1^{\otimes a_1}\otimes\cdots\otimes L_k^{\otimes a_k}$ the line bundle on $F$, where $a_i-a_{i+1}\geq l$ for $1\leq i\leq k-1$ and $L|_X$ its restriction on $X$.

We  compare the evaluation map of $L $ on $X$ and $F=\mathrm{Fl}(n_1,\ldots, n_k; V)$. Let $M_L$ be the kernel of  $$ev_L: H^0(F, L )\otimes \cO_{F}\rightarrow L $$ and let $M_{L|X}$ be the kernel of $$ev_{L|_X}: H^0(X, L|_X)\otimes \cO_{X}\rightarrow L|_X.$$ By Lemma \ref{restriction}, since the restriction map is equivariant under the action of the underlying symplectic group which is reductive, we may regard $H^0(F, L)$ as a direct sum of $H^0(X, L|_X)$ and $H^0(F, L\otimes \cI_X)$, where $\cI_X$ is the ideal sheaf on $F$ defining $X$.  We now restrict the evaluation map $ev_L$ on $X$ and get $$0\rightarrow M_L|_X\rightarrow H^0(F, L)\otimes \cO_X\xrightarrow{ev_L|_X} L|_X\rightarrow 0.$$ By the above decomposition and since the evaluation maps $ev_L|_X$ and $ev_{L|_X}$ are equivariant under the symplectic group, we see that $ev_L|_X$ splits as the direct sum of $ev_{L|_X}$ with the zero map $H^0(F, L\otimes \cI_X)\otimes \cO_X\rightarrow 0$. Thus
$M_L|_{X}\simeq M_{L|_X}\bigoplus (H^0(F, L\otimes \cI_X)\otimes \cO_X).$ Hence, in order to show that $L|_{X}$ satisfies Property $(N_p)$, it suffices to prove that $$H^1(X, \wedge^i(M_L|_{X})\otimes L^{\otimes h})=0,$$ for $1\leq i\leq p+1$ and $h\geq 1$.
Combining the sequence (\ref{resolution}), we only need to show
\begin{eqnarray}\label{vanishing}H^{1+j}(F, \wedge^iM_L\otimes L^{\otimes h}\otimes \wedge^jW)=0,\end{eqnarray} for $1\leq i\leq p+1$, $h\geq 1$ and  $j\geq 0$.

By (\ref{quotients}), we just need to show that
$$H^{1+j}(F, \wedge^iM_L\otimes L^{\otimes h}\otimes S^{\rho_2}\cQ_2\otimes\cdots\otimes S^{\rho_{k+1}}\cQ_{k+1})=0,$$ for  $1\leq i\leq p+1$, $h\geq 1$, $j\geq 0$, and $\rho_2, \ldots, \rho_{k+1}$ as in (\ref{quotients}).

Note that $\wedge^iM_L$ is a direct  summand of $\otimes^iM_L$ and  $M_L$ admits a filtration whose associated graded bundles are $N_i$ for $1\leq i\leq k$, which are described in Subsection 2.3. Thus it suffices to show that
$$H^{1+j}(F, N_{t_1}\otimes\cdots\otimes N_{t_i}\otimes L^{\otimes h}\otimes S^{\rho_2}\cQ_2\otimes\cdots\otimes S^{\rho_{k+1}}\cQ_{k+1})=0$$ for $1\leq t_1,\ldots, t_i\leq k$.

By Lemma \ref{weyman}, the above vanishing is a consequence of the vanishing
\begin{eqnarray}\label{vanishing-resolution}\;\;\;\;\;\;\;\;\;H^{1+j+(j_1\cdots+j_i)-i}(F, C(t_1)_{j_1}\otimes\cdots\otimes C(t_i)_{j_i}\otimes L^{\otimes h}\otimes S^{\rho_2}\cQ_2\otimes\cdots\otimes S^{\rho_{k+1}}\cQ_{k+1})=0,
\end{eqnarray}
for $1\leq i\leq p+1$, $j\geq 0$, and $ j_1,\ldots, j_i > 0$. We may let $h=1$ in the following, since the calculation goes the same way and $h=1$ is indeed the  most delicate case.

Recall that $C(t_s)_{j_s}$ is a direct sum of vector bundles which can be expressed as
$$S^{\mu}\cQ_{t_s+1}\otimes S^{\nu}(\Sigma_0/\Sigma_{t_s+1})\otimes \tilde{L}_{t_s+1},$$ where $|\mu|=j_s$. Thus by Lemma \ref{finalvanishing},  in order to prove   (\ref{vanishing-resolution}),
it suffices to show that
given any partitions $\mu_s$ of weight $j_s$ for $1\leq s\leq i$ and a direct summand $$S^{\alpha_2}\cQ_2\otimes \cdots\otimes S^{\alpha_{k+1}}\cQ_{k+1}$$  of
$$\big(S^{\mu_1}\cQ_{t_1+1}\otimes\cdots\otimes S^{\mu_i}\cQ_{t_i+1}\big)\otimes S^{\rho_2}\cQ_2\otimes\cdots\otimes S^{\rho_{k+1}}\cQ_{k+1},$$ we have  $$1+j+(j_1+\cdots+j_i)-i> \sum_{2\leq p\leq k+1}\sum_{1\leq i\leq s_j}a_{pi}-ls-\sum_{2\leq j\leq k+1}s_j^2,$$  for any integers $0\leq s_2\leq r_2, \ldots, 0\leq s_{k+1}\leq r_{k+1}$ with $s:=s_2+\ldots+s_{k+1}>0$, where $a_{pi}$ is the $i$-th entry of $\alpha_p$.


 By the Pieri rule, we have that, \begin{eqnarray*} \sum_{2\leq p\leq k+1}\sum_{1\leq i\leq s_j}a_{pi}& \leq  &(|\mu_1|+\cdots +|\mu_i|)+\sum_{2\leq p\leq k+1}\sum_{1\leq i\leq s_j}b_{pi}\\
&=& (j_1+\cdots +j_i)+\sum_{2\leq p\leq k+1}\sum_{1\leq i\leq s_j}b_{pi},
\end{eqnarray*} where $b_{pi}$ is the $i$-th entry of $\rho_p$.   By Lemma \ref{precise1}, we have
\begin{eqnarray}\label{difference} \sum_{2\leq p\leq k+1}\sum_{1\leq i\leq s_j}b_{pi} \leq j+\frac{s(s-1)}{2}.
\end{eqnarray}
Combining the above inequalities, it rests to prove the inequality
$$1+j+(j_1+\cdots+j_i)-i> (j_1+\cdots+j_i)+j+\frac{s(s-1)}{2}-ls-\sum_{2\leq j\leq k+1}s_j^2.$$ Since $1\leq i\leq p+1$, we just need that
\begin{eqnarray}l\geq \frac{p+1}{s}+\frac{s-1}{2}-\frac{\sum_{2\leq j\leq k+1}s_j^2}{s},
\end{eqnarray}
for any any integers $0\leq s_2\leq r_2, \ldots, 0\leq s_{k+1}\leq r_{k+1}$ with $s:=s_2+\ldots+s_{k+1}>0$.

We can now finish the proof of Theorem \ref{symplectic}.

When $k\leq 2$, $\frac{s}{2}\leq \frac{\sum_{2\leq j\leq k+1}s_j^2}{s}$. Thus, as soon as $l\geq p$, $L$ satisfies Property $(N_p)$.

In general, $1\leq  \frac{\sum_{2\leq j\leq k+1}s_j^2}{s}$, thus as soon as $l\geq \frac{p+1}{s}+\frac{s-3}{2}$, $L$ satisfies property $(N_p)$. Since $$1\leq s\leq \sum_{2\leq j\leq k+1}r_j=n_1$$ and the function
\begin{eqnarray*}&\mathbb R&\rightarrow \mathbb R\\
&s&\rightarrow \frac{p+1}{s}+\frac{s-3}{2}
\end{eqnarray*}
decreases for $s\in (1, \sqrt{p+1})$ and increases for $s\in (\sqrt{p+1}, +\infty)$,
  we conclude that when $l\geq \mathrm{max}\{p, \frac{p+1}{n_1}+\frac{n_1-3}{2}\}$, $L$ satisfies property $(N_p)$. We also note that when $p\geq \frac{n_1}{2}-1$, $p\geq \frac{p+1}{n_1}+\frac{n_1-3}{2}$.
  \end{proof}

  \begin{rema}
  If $X=\mathrm{SFl}(1,2,\ldots, n; V)$ is the full symplectic variety, the estimation in the last paragraph of the proof of Theorem \ref{symplectic} seems to be optimal.
  However, we can sometimes refine the estimation on $l$ in other cases.
  The proof shows that $L$ on $X$ satisfies Property $(N_p)$ if
  $$l\geq \max\{\frac{p+1}{s}+\frac{s-1}{2}-\frac{\sum_{2\leq j\leq k+1}s_j^2}{s}\mid 0\leq s_i\leq r_i\;\mathrm{and}\; 1\leq s\leq n_1\}.$$

   We now  consider $X=\mathrm{SFl}(3, 5, 6; V)$, where $V$ is of dimension $12$.
   The Picard number of $X$ is $3$ and $\dim X=47$. We also have $r_2=1$, $r_3=2$ and $r_4=3$.
   Then an elementary calculation shows that $L$ satisfies Property $(N_p)$ whenever $l\geq p$.    \end{rema}
  \subsection{Orthogonal flag varieties}

  In the type $B_n$ and $D_n$ cases, a crucial difference is that  the Picard group $X$ and the Picard group of $F$ could be slightly different and our arguments work only for  ample line bundles coming from $F$.

   \begin{theo}\label{quadric}
   Assume that $X$  is a rational homogeneous variety of type $B_n$ or $D_n$ and let $\iota: X\hookrightarrow F$ be the embedding of $X$ to the corresponding flag variety explained in Section 3.
   Assume that  $L$ can be expressed as  $\iota^*(H^{\otimes l}\otimes M)$ where $H$ and $M$ are respectively  ample and nef  line bundles on $F$. Then $L$ satisfies Property $(N_p)$, if one of the following conditions hold:
\begin{itemize}
\item[(1)] $k\leq 2$ and $l\geq p+1$;
\item[(2)]  $l\geq p+1\geq \frac{n_1}{2}-1$;
\item[(3)]  $l\geq \mathrm{max}\{p+1, \frac{p+1}{n_1}+\frac{n_1-1}{2}\}$.
\end{itemize}
 \end{theo}

 The proof is identical to the proof of Theorem \ref{symplectic}. We just need to replace the inequality (\ref{difference}) by
 $$ \sum_{2\leq p\leq k+1}\sum_{1\leq i\leq s_j}b_{pi} \leq j+\frac{s(s+1)}{2}.$$

We remark that when $n_1<n$ in the $B_n$ case or $n_1<n-1$ in the $D_n$ case, the restriction map $\Pic(F)\rightarrow \Pic(X)$ is bijective (see Subsection 3.2 and 3.3) and thus we may remove $"\iota^*"$ in the assumption of Theorem \ref{quadric}.
 \subsection{$G_2$-flag varieties}

  We apply the same strategy to study the syzygy of homogeneous varieties of type $G_2$. We need the commutative diagram in Section 3:
 \begin{eqnarray*}
 \xymatrix{
X \ar@{^{(}->}[d]& \mathbb P_X(\Sigma)\ar[l]\ar@{^{(}->}[d]\ar[r] & Q\ar@{=}[d]\\
\mathrm{OGr}(2, V)\ar@{^{(}->}[d]&\mathrm{OFl}(1, 2; V)\ar[l]\ar@{^{(}->}[d]\ar[r] & Q\ar@{^{(}->}[d]\\
\mathrm{Gr}(2, V)&\mathrm{Fl}(1, 2; V)\ar[l]^p\ar[r]& \mathbb P(V)=\mathrm P^6,}
 \end{eqnarray*}
 where the varieties on the top line are rational homogeneous varieties of type $G_2$. Since $Q$ is a $5$-dimensional quadric hypersurface, we just need to deal with $X$ and $\mathbb P_X(\Sigma)$.

Let $\Sigma$ and $\cQ$ be the tautological sub-bundle and quotient bundle on $\mathrm{Gr}(2, V)$. Let $L_G$ be the Pl\"ucker polarization on $\mathrm{Gr}(2, V)$.  We know that $X$ is defined as the zero locus of a general section of the vector bundle $\cQ^*\otimes L_G$ on $\mathrm{Gr}(2, V)$. Then $-K_X=(L_G|_X)^{\otimes 3}$ and $\Pic(X)$ is generated by $L_G|_X$.
 Similarly, $\mathbb P_X(\Sigma)$ is defined as the zero locus of a general section of the vector bundle $p^*(\cQ^*\otimes L_G)$ on $\mathrm{Fl}(1, 2; V)$ and
the restriction map $$\Pic(\mathrm{Fl}(1, 2; V))\rightarrow \Pic( \mathbb P_X(\Sigma))$$ is also bijective.
 \begin{lem}\label{G2lemma} For any ample line bundle $L$ on $\mathrm{Gr}(2, V)$ (resp. $\mathrm{Fl}(1, 2; V)$), the restriction map $$H^0(\mathrm{Gr}(2, V), L)\rightarrow H^0(X, L|_X)$$ (resp. $H^0(\mathrm{Fl}(1, 2; V), L)\rightarrow H^0(\mathbb P_X(\Sigma), L|_{\mathbb P_X(\Sigma)}))$ is surjective.
 \end{lem}
 \begin{proof}
 We may argue as in the proof of Lemma \ref{restriction}. It is also easy to write down a direct proof via the Borel-Weil-Bott theorem.

 \end{proof}
 \begin{theo}\label{G2}
 Let $L$ be a an ample line bundle on $X$ or $ \mathbb P_X(\Sigma)$. Assume that $L=H^{\otimes l}\otimes M$, where $H$ is ample, $l\geq 1$, and $M$ is nef. Then $L$ satisfies Property $(N_p)$ when $l\geq p$.
 \end{theo}
 \begin{proof}
We apply the same strategy of the proof of Theorem \ref{symplectic}. By Lemma \ref{G2lemma}, we may assume that $p\geq 1$. We first consider $X\hookrightarrow \mathrm{Gr}(2, V)$. We also regard $L$ as an ample line bundle on $ \mathrm{Gr}(2, V)$ and let $M_L$ be its kernel bundle on $ \mathrm{Gr}(2, V)$. Since the Picard number of a Grassmannian is $1$, we may assume that $L=L_G^{\otimes l}$.

We then need to show that
$$H^{1+j}( \mathrm{Gr}(2, V), \wedge^iM_L\otimes \wedge^j(\cQ\otimes L_G^{-1})\otimes L)=0.$$ We apply Lemma \ref{weyman} and the Borel-Weil-Bott theorem. Note that on $\mathrm{Gr}(2, V)$,  Lemma \ref{weyman} simply says that $M_L$ admits a resolution
$$\cdots\rightarrow C(j)\rightarrow \cdots\rightarrow C(1)\rightarrow M_L\rightarrow 0,$$
where $C(j)$ is a direct sum of copies of  $S^{\rho}\Sigma$ with $|\rho|=j$. Thus we are reduced to show that
\begin{eqnarray}\label{newvanishing}
H^{1+j-i+\sum_{1\leq s\leq i}|\rho_s|}(S^{\rho_1}\Sigma\otimes\cdots\otimes S^{\rho_i}\Sigma\otimes  \wedge^j(\cQ\otimes L_G^{-1})\otimes L)=0,
\end{eqnarray}
 where $0\leq j\leq 5$, $1\leq i\leq p+1$, and $|\rho_s|\geq 1$ for $1\leq s\leq i$.

By the Borel-Weil-Bott Theorem, given a partition $(a_1, a_2)$ and an integer $m>0$,    $$H^m(S^{(a_1, a_2)}\Sigma\otimes \wedge^j(\cQ\otimes L_G^{-1})\otimes L)$$ is nonzero only in the following cases:
\begin{itemize}
\item[1)] $m=5$, $a_1-6>l-j$, and $a_2-7<l-j-5$;
\item[2)] $m=10$, $a_1-6>a_2-7>l-j$;
\item[3)]$m=5-j$ when $1\leq j\leq 4$, $a_1-6=l-2j$, and $a_2-7<l-j-5$;
\item[4)] $m=10-j$ when $1\leq j\leq 4$, $a_1-6>l-j$, and $a_2-7=l-2j$.

\end{itemize}
 Let  $S^{(a_1, a_2)}\Sigma$ be a direct summand of $S^{\rho_1}\Sigma\otimes\cdots\otimes S^{\rho_i}\Sigma$. Note that $a_1+a_2=\sum_{1\leq k\leq i}|\rho_k|\geq i$ and $l\geq p\geq i-1$.
 Assume to the contrary that $$H^{1+j-i+a_1+a_2}(S^{(a_1, a_2)}\Sigma\otimes \wedge^j(\cQ\otimes L_G^{-1})\otimes L )\neq 0,$$ we now deduce a contradiction.

  If $1+j-i+a_1+a_2=5$,  we have by $1)$ that
  $$5=1+j-i+a_1+a_2>a_2+1+j-i+(l+6-j)) \geq a_2+6>5,$$ which is impossible.

   If $1+j-i+a_1+a_2=10$,  we have by $2)$ that
   $$10=1+j-i+a_1+a_2\geq 1+j-i+2(l+8-j)=17-j+2l-i \geq l+16-j\geq l+11,$$ which is impossible.

   If $1+j-i+a_1+a_2=5-j$ for some $1\leq j\leq 4$, we have by $3)$ that
   $$5-j=1+j-i+a_1+a_2=1+j-i+(l+6-2j)+a_2=l+7-i-j+a_2\geq 6-j+a_2,$$ which is impossible.

   Finally, if $1+j-i+a_1+a_2=10-j$  for some $1\leq j\leq 4$, we have by $4)$ that
   $$10-j=1+j-i+a_1+a_2>1+j-i+(l+6-j)+(l+7-2j)=14+2l-i-2j,$$ which implies that $j>3+l\geq 4$. This is again absurd.\\\\

We now consider $ \mathbb P_X(\Sigma)$. Let $L=L_1^{\otimes a}\otimes L_2^{\otimes b}$ with $a-b\geq l$ and $b\geq l$.  We recall the notations and  results in Subsection 2.3.  In this case, since the ambient variety is $\mathrm{Fl}(1, 2; V)$, we have $\Sigma_3=0$, $\cQ_3=\Sigma_2=(L_1\otimes L_2)^{-1}$, $\cQ_2=L_2$, and $\cQ_1=p^*\cQ$. We have $H^0(L)\simeq S^{\alpha}V$, where $$\alpha=(\overbrace{a,\ldots, a}^5, b, 0).$$

Moreover, $M_L$ sits in a short exact sequence $$0\rightarrow N_2\rightarrow M_L\rightarrow N_1\rightarrow 0.$$ By Lemma \ref{weyman}, $N_2$ admits a resolution $$\cdots\rightarrow C(2)_j\rightarrow \cdots\rightarrow C(2)_1\rightarrow N_2\rightarrow 0,$$ where $C(2)_j$ is a direct sum of copies of  $(L_1\otimes L_2)^{\otimes (-j)}$ since $\cQ_3=(L_1\otimes L_2)^{-1}$ and $\Sigma_3=0$. By the definition of the Littlewood-Richardson coefficient in Lemma \ref{weyman}, we know by the Pieri rule that $1\leq j\leq a$. Similarly,  $N_1$ admits a resolution $$\cdots\rightarrow C(1)_j\rightarrow \cdots\rightarrow C(1)_1\rightarrow N_1\rightarrow 0,$$ where, since $\alpha_{\geq 1}=(a-b, a-b, a-b, a-b, a-b, 0)$, we know by the Pieri rule that $1\leq j\leq a-b$ and $$C(1)_j=L_2^{\otimes j}\otimes S^{(a-b, a-b, a-b, a-b, a-b-j, 0)}(\Sigma_0/\Sigma_2)\otimes (L_1\otimes L_2)^{\otimes b}.$$   We need to show that
\begin{multline}\label{11}H^{1+j-i+\sum_{1\leq l\leq i}j_l}(\mathrm{Fl}(1, 2; V), C(t_1)_{j_1}\otimes\cdots \otimes C(t_l)_{j_l}\otimes\cdots\otimes C(t_i)_{j_i}\otimes\\  \wedge^jp^*(\cQ\otimes L_G^{-1})\otimes L)=0\end{multline} for $0\leq j\leq 5$, $1\leq i\leq p+1$, $t_1,\ldots, t_i=1$ or $2$, and $j_1,\ldots, j_i\geq 1$.

We denote $J=\sum_{1\leq l\leq i}j_l$,  $V=\{l\mid 1\leq l \leq i \; \text{and}\; t_l=1\}$, and $J_1=\sum_{l\in V}j_l$.

Let $S^{\rho}(\Sigma_0/\Sigma_2)\otimes L_1^{\otimes s}\otimes L_2^{\otimes t}$ be a direct summand of $C(t_1)_{j_1}\otimes\cdots\otimes C(t_i)_{j_i}$. From the discussions above about $C(t_l)_{j_l}$, we see that
 \begin{eqnarray}\label{s-t} && s=b|V|-J+J_1\\\nonumber && t=b|V|-J+2J_1.
 \end{eqnarray}
We observe that $s-t=-J_1\geq -J$.

Therefore, in order to prove (\ref{11}), we are reduced to show the vanishing that
\begin{eqnarray}\label{12}H^{1+j-i+J}(\wedge^j\cQ_1\otimes L_1^{\otimes (a-j+s)}\otimes L_2^{\otimes (b+t)}\otimes S^{\rho}(\Sigma_0/\Sigma_2))=0.\end{eqnarray}


We now apply the Borel-Weil-Bott Theorem and the same induction argument in the proof of  Proposition \ref{finalvanishing}.

For each weight $\rho$ of positive weight, $S^{\rho}(\Sigma_0/\Sigma_2)$ is a direct summand of $(\Sigma_0/\Sigma_2)\otimes S^{\rho'}(\Sigma_0/\Sigma_2)$ where $|\rho'|=|\rho|-1$. Since $\Sigma_0$ is a trivial vector bundle and $\Sigma_2=(L_1\otimes L_2)^{-1}$, by induction on $|\rho|$, it suffices to show that
\begin{eqnarray*}H^{1+j-i+J+K}(\wedge^j\cQ_1\otimes L_1^{\otimes (a-j+s)}\otimes L_2^{\otimes (b+t)}\otimes (L_1\otimes L_2)^{\otimes (-K)})=0\end{eqnarray*} for each $K\geq 0$, to prove (\ref{12}).

By the Borel-Weil-Bott theorem, it is reduced  to study the sequence
\begin{eqnarray*}
&&(b_1, b_2, b_3, b_4, b_5, b_6, b_7):=\\&&(\overbrace{a+s+1-j-K,\ldots, a+s+1-j-K}^j, \overbrace{a+s-K-j, \ldots, a+s-K-j}^{5-j}, b+t-K, 0)-\\&&(1,2,3,4,5,6,7)=\\
&&(\overbrace{a+s-j-K,\ldots, a+s-K-2j+1}^j, \overbrace{a+s-2j-K-1, \ldots, a+s-j-K-5}^{5-j},\\&& b+t-K-6, -7)\end{eqnarray*} of length  $7$.

If the above sequence has two identical entries or is strictly decreasing, nothing needs to be proved. Otherwise, the number of inversions of the sequence is nonzero, then it can only be $1$, $5$, $6$, $10$, $11$ when $j=0$ or $5$ and some more numbers $5-j$, $6-j$, $10-j$ and $11-j$ when $1\leq j\leq 4$. We finally show that in  all these possibilities, $1+j-i+J+K$ is not the number of inversions of the above sequence. This is an elementary calculation but we still write down the details for the completeness of the proof.

If $1+j-i+J+K=1$ is the number of inversions and the set of inversions is just $\{b_6<b_7\}$, we have $K=j=0$ and $J=i$. Moreover, $b+t-6<-7 $. Recall (\ref{s-t}) that $t=b|V|-J+2J_1\geq -J=-i\geq -p-1$. Thus $p\leq l\leq b<-t-1\leq p$, which is a contradiction.

If $1+j-i+J+K=5$ is the number of inversions and the set of inversions is $\{b_6>b_e\mid  1\leq e\leq 5\}$, $J=4+i-j-K\leq 5+p-j$. In this case, we also have $b+t-K-6\geq a+s-K-j$. Then $J\geq t-s\geq (a-b)+6-j\geq p+6-j$, which is a contradiction.

If $1+j-i+J+K=6$ is the number of inversions and the set of inversions is $\{b_7>b_e\mid  1\leq e\leq 6\}$, $K=5-j+i-J\leq 6+p-j-J$. On the other hand, we have $-7\geq a+s-K-j$. Thus $K\geq a+s+7\geq 2p-j-J+7$, which is a contradiction.

If $1+j-i+J+K=10$ is the number of inversions and the set of inversions is $\{b_6>b_e, b_7>b_e\mid 1\leq e\leq 5\}$, we have $K=9-j+i-J\leq 10+p-j-J$. In this case, we also have $-7\geq a+s-K-j$ and $b+t-6\geq a+s-j+1$. Thus $$10+p-j-J\geq K\geq a+s+7-j\geq 2p+7+b|V|-J+J_1-j$$ and we have $J_1\leq 3-p$. On the other hand, from the inequality $b+t-6\geq a+s-j+1$, we get $p+2\leq p+7-j\leq J_1$, which is a contradiction.

If $1+j-i+J+K=11$ is the number of inversions and the set of inversions is $\{b_6>b_e, b_7>b_e, b_7>b_6\mid 1\leq e\leq 5\}$, $K=10-j+i-J\leq 11+p-j-J$. Moreover, since $-7>a+s-j-K$, we have $K\geq a+s-j+8\geq 2p+J_1-J-j+8$. Thus $J_1\leq 3-p$.
On the other hand, $b+t-6\geq a+s-j$ and thus $J_1\geq p+6-j$. We thus have a contradiction when $j\leq 4$.  When $j=5$, we can improve the above estimation. Indeed, we have $b+t-6\geq a+s-4$ and thus $J_1\geq p+2$ and we still have a contradiction.
\\

We then assume that $1\leq j\leq 4$.

If $1+j-i+J+K=5-j$ is the number of inversions and the set of inversions is $\{b_6>b_e\mid  j+1\leq e\leq 5\}$, $2j=4+i-J-K\leq 5+p-J$.  We also have that $b+t-K-6=a+s-K-2j$ in this case. Then $2j=6+(a-b)+(s-t)\geq l+6-J\geq p+6-J$, which is a contradiction.

If $1+j-i+J+K=6-j$ is the number of inversions and the set of inversions is $\{b_7>b_e\mid  j+1\leq e\leq 6\}$, $2j=5+i-J-K\leq 6+p-J-K$. In this case, we have $-7=a+s-K-2j$. Thus $2j=a+s+7-K\geq 2p-J+7-K$, which is a contradiction.

If $1+j-i+J+K=10-j$ is the number of inversions and the set of inversions is $\{b_6>b_e, b_7>b_f\mid 1\leq e\leq 5\; \text{and} \;j+1\leq f\leq 5\}$, $2j=9+i-J-K\leq 10+p-J-K$. Moreover, we have $-7=a+s-K-2j$ and $b+t-6\geq a+s-j+1$. Thus $J_1=t-s\geq a-b+7\geq p+7$ and $2j=a+7+s-K\geq 2p+7+J_1-J-K\geq 3p+14-J-K$, which is a contradiction.

If $1+j-i+J+K=11-j$ is the number of inversions and the set of inversions is $\{b_7>b_e, b_6>b_f\mid 1\leq e\leq 5\; \text{and} \;j+1\leq f\leq 5\}$, we have $2j=10+i-J-K\leq 11+p-J-K$. Moreover, $-7\geq a+s-K+1$ and $b+t-6=a+s-2j$. Thus $-7\geq 2p+J_1-J-K+1$ and $J+K\geq 2p+J_1+8$. On the other hand, $2j=(a-b)+(s-t)+6\geq p+6-J_1$. Thus $p+6-J_1\leq 2j\leq 11+p-J-K\leq 11+p-2p-J_1-8=3-p-J_1$, which is a contradiction.

 \end{proof}
 
\end{document}